\documentclass{amsart}

\newtheorem{theorem}{Theorem}
\newtheorem{lemma}[theorem]{Lemma}
\newtheorem{proposition}[theorem]{Proposition}

\theoremstyle{definition}

\usepackage{amscd,amssymb}

\begin{document}

\title[Trace cocharacters of two $4\times 4$ matrices]
{Multiplicities in the trace cocharacter sequence of
two $4\times 4$ matrices}

\author[Vesselin Drensky and Georgi Genov]
{Vesselin Drensky and Georgi K. Genov}
\address{Institute of Mathematics and Informatics,
Bulgarian Academy of Sciences, 1113 Sofia, Bulgaria}
\email{drensky@math.bas.bg, gguenov@hotmail.com}

\thanks{Partially supported by Grant
MM-1106/2001 of the Bulgarian Foundation for Scientific Research.}

\begin{abstract}
We find explicitly the generating functions of the
multiplicities in the pure and mixed trace
cocharacter sequences of two $4\times 4$ matrices over a field of
characteristic 0. We determine the asymptotic behavior of the multiplicities and
show that they behave as polynomials of 14th degree.
\end{abstract}
\maketitle

\section*{Introduction}

Let us fix an arbitrary field $F$ of characteristic 0 and two integers
$n,d\geq 2$. Consider the $d$ generic $n\times n$ matrices
$X_1,\ldots,X_d$. We denote by $C$ the pure trace algebra
generated by the traces of all products
$\text{\rm tr}(X_{i_1}\cdots X_{i_k})$, and by $T$ the
mixed trace algebra generated by $X_1,\ldots,X_d$
and $C$, regarding the elements of $C$ as scalar matrices.
The algebra $C$ coincides with the algebra of
invariants of the general linear group $GL_n(F)$ acting by
simultaneous conjugation on $d$ matrices of size $n$. The algebra
$T$ is the algebra of matrix concominants under a suitable action of
$GL_n(F)$. See e.g. the books \cite{J}, \cite{F3}, or \cite{DF} as
a background on $C$ and $T$ and their numerous applications.

The algebras $C$ and $T$ are graded by multidegree.
The Hilbert (or Poincar\'e) series of $C$ is
\[
H(C)=H(C,t_1,\ldots,t_d)=\sum\dim C^{(k)}t_1^{k_1}\cdots
t_d^{k_d},
\]
where $C^{(k)}$ is the homogeneous component of multidegree
$k=(k_1,\ldots,k_d)$. In the same way one defines the Hilbert series of
$T$. These series are symmetric functions and
decompose as infinite linear combinations of Schur functions
$S_{\lambda}(t_1,\ldots,t_d)$,
\[
H(C)=\sum m_{\lambda}(C)S_{\lambda}(t_1,\ldots,t_d),
\]
where $\lambda=(\lambda_1,\ldots,\lambda_d)$ is a partition in not more than
$d$ parts, and the $m_{\lambda}(C)$s are nonnegative integers which are 0 if
$d>n^2$ and $\lambda_{n^2+1}>0$.
A similar expression holds for the Hilbert series of $T$. The
multiplicities $m_{\lambda}(C)$ and $m_{\lambda}(T)$ have important
combinatorial properties and ring theoretical meanings. In particular,
they are equal to the multiplicities of the irreducible $S_k$-characters
in the sequences of pure and mixed trace cocharacters, respectively,
and give estimates for the multiplicities in the ``ordinary'' cocharacters
of the polynomial identities of the $n\times n$ matrix algebra.

The multiplicities of $m_{\lambda}(C)$ and $m_{\lambda}(T)$ are
explicitly known in very few cases. For $n=2$ and any $d$, Formanek \cite{F1}
showed that
\[
m_{\lambda}(T)=(\lambda_1-\lambda_2+1)(\lambda_2-\lambda_3+1)(\lambda_3-\lambda_4+1),
\]
if $\lambda=(\lambda_1,\lambda_2,\lambda_3,\lambda_4)$, and
$m_{\lambda}(T)=0$ if $\lambda_5\not=0$. He also found a more complicated
expression for the multiplicities in the Hilbert series of $C$.
With some additional work, based on representation theory of symmetric and general linear
groups, these results can be also derived from the paper by Procesi \cite{P}.
The Hilbert series of $C$ and $T$ can be expressed as multiple integrals,
evaluated by Teranishi \cite{T1, T2} for $C$ and $n=3,4$ and $d=2$. Van
den Bergh \cite{V} suggested graph theoretical methods for
calculation of $H(C)$ and $H(T)$. Berele and Stembridge
\cite{BS} calculated the Hilbert series of $C$ and $T$ for
$n=3$, $d\leq 3$ and for $n=4$, $d=2$, correcting also some
typographical errors in the expression of $H(C)$ for $n=4$ and $d=2$ in
\cite{T2}. Using the Hilbert series of $C$, $n=3$, $d=2$, Berele
\cite{B1} found an asymptotic expression of
$m_{(\lambda_1,\lambda_2)}(C)$. The explicit form of the
generating function of $m_{(\lambda_1,\lambda_2)}(C)$ was found by
the authors of the present paper \cite{DG1}
correcting also a technical error in
\cite{B1}. Later they \cite{DG2} suggested a method to find the
coefficients of the Schur functions in the expansion of a class of
rational symmetric functions in two variables. Jointly with Valenti \cite{DGV}
they have used the expression of the Hilbert series $H(T)$  found
by Berele and Stembridge \cite{BS} and have calculated explicitly the multiplicities
$m_{(\lambda_1,\lambda_2)}(T)$ for $n=3$ and $d=2$.

The purpose of the present paper is to calculate, for $n=4$ and $d=2$, the
generating functions of $m_{(\lambda_1,\lambda_2)}(C)$ and
$m_{(\lambda_1,\lambda_2)}(T)$. In principle, this allows to give explicit expressions for
the multiplicities. Since the formulas are quite complicated,
we prefer to give the asymptotics only, as in \cite{B1}. It turns out that the
multiplicities $m_{(\lambda_1,\lambda_2)}(C)$ and
$m_{(\lambda_1,\lambda_2)}(T)$ behave
as polynomials of degree 14 in $\lambda_1$ and $\lambda_2$. As in \cite{DGV},
our approach is to apply the methods of \cite{DG2} to
the explicit form of the Hilbert series of $C$ and $T$ found in \cite{BS}.
As in \cite{DGV}, results of
Formanek \cite{F1, F2} imply that the values of the
multiplicities $m_{\lambda}(M_4(F))$ for
$\lambda=(\lambda_1,\ldots,\lambda_{16})$ coincide with these of
$m_{(\lambda_1-\lambda_{16},\lambda_2-\lambda_{16})}(T)$ when
$\lambda_3=\cdots=\lambda_{16}\geq 2$.

Carbonara, Carini and Remmel \cite{CCR} determined, for any $n$,
the behaviour of $m_{\lambda}(C)$, when $\lambda_2,\ldots,\lambda_{n^2}$ is fixed and
$\lambda_1$ sufficiently large. Our results give more detail, for $n=4$, in this direction.
A general result of Berele \cite{B2} describes the multiplicities of Schur functions for
the class of rational symmetric functions with denominators which are products of
binomials $1-t_1^{k_1}\cdots t_d^{k_d}$, in any number of variables. It covers the Hilbert series
of relatively free algebras and pure and mixed trace algebras. Again, our results give some detail
in the partial case which we consider in this paper.

\section{Preliminaries}

We refer to the book by Macdonald \cite{M} for a background on theory of symmetric
functions. We shortly summarize the facts we need for our exposition.
We consider the algebra $\text{\rm Sym}[[x,y]]$
of formal power series which are symmetric functions in two variables over $\mathbb C$.
Every element $f(x,y)\in \text{\rm Sym}[[x,y]]$ can be expressed in a unique way
as an infinite linear combination
\[
f(x,y)=\sum_{\lambda}m(\lambda)S_{\lambda}(x,y),
\]
where $S_{\lambda}(x,y)$ is the Schur function related with the partition
$\lambda=(\lambda_1,\lambda_2)$,
and $m(\lambda)\in \mathbb C$ is the multiplicity of $S_{\lambda}(x,y)$
in the decomposition
of $f(x,y)$. The Schur function $S_{\lambda}(x,y)$ has the form
\[
S_{\lambda}(x,y)=(xy)^{\lambda_2}(x^p+x^{p-1}y+\cdots+xy^{p-1}+y^p)=
\frac{(xy)^{\lambda_2}\left(x^{p+1}-y^{p+1}\right)}{x-y},
\]
where we have denoted $p=\lambda_1-\lambda_2$.
In \cite{DG1} we introduced the multiplicity series of
$f(x,y)$
\[
M(f)(t,u)=\sum_{\lambda_1\geq\lambda_2\geq 0}
m(\lambda_1,\lambda_2)t^{\lambda_1}u^{\lambda_2} \in {\mathbb C}[[t,u]].
\]
Introducing a new variable $v=tu$, it accepts the more convenient form
\[
M'(f)(t,v)=M(f)(t,u)=\sum_{\lambda_1\geq\lambda_2\geq 0}
m(\lambda_1,\lambda_2)t^{\lambda_1-\lambda_2}v^{\lambda_2} \in
{\mathbb C}[[t,v]].
\]
The mapping $M':\text{\rm Sym}[[x,y]]\to {\mathbb C}[[t,v]]$ is a continuous
linear bijection.
The relation between symmetric functions and their multiplicity
series is given, see (11) in \cite{DG1}, by
\[
f(x,y)=\frac{xM'(f)(x,xy)-yM'(f)(y,xy)}{x-y}. \eqno(1)
\]
The following theorem
from \cite{DG2} describes the symmetric functions with rational
multiplicity series and gives hints how to calculate the multiplicity
series for such functions.

\begin{theorem}
The multiplicity series of $f(x,y)\in \text{\rm Sym}[[x,y]]$ is rational
if and only if $f(x,y)$ has the form
\[
f(x,y)=\frac{p(x,xy)+p(y,xy)}{q(x,xy)q(y,xy)},
\]
where $p(x,z),q(x,z)$ are polynomials in $x$ with coefficients
which are rational functions in $z$. Then
\[
M_1(f;t,v)=\frac{h(t,v)}{q(t,v)},
\]
where $h(t,v)\in {\mathbb C}(v)[t]$ and
$\text{\rm deg}_th\leq\max(\text{\rm deg}_xp,\text{\rm deg}_xq-2)$.
\end{theorem}
The polynomial $h(t,v)$ can be found from (1)
by the method of unknown coefficients.
The calculations can be simplified if we know the decomposition of $q(x,z)$.
We need the following easy lemma which is a slight generalization
of Lemma 14 in \cite{DG1}.

\begin{lemma}
Let $K$ be any field, let $\xi$ be an arbitrary element from $K$
and let $f(w), g(w)\in K[w]$ be two polynomials such that
$f(1/\xi),g(1/\xi)\not=0$. Then, in the decomposition as a sum of elementary fractions,
\[
\frac{f(w)}{(1-\xi w)^kg(w)}=
\frac{\alpha_k}{(1-\xi w)^k}+\frac{\alpha_{k-1}}{(1-\xi w)^{k-1}}
+\cdots+\frac{\alpha_1}{1-\xi w}+
\frac{b(w)}{g(w)}+c(w),
\eqno(2)
\]
where $\alpha_1,\ldots,\alpha_k\in K$ and $b(w),c(w)\in K[w]$, the coefficient $\alpha_k$
has the form
\[
\alpha_k=\frac{f(1/\xi)}{g(1/\xi)}.
\]
\end{lemma}

\section{Main Results}

In the sequel we fix $n=4$ and $d=2$ and denote, respectively, by
$C$ and $T$ the pure and mixed trace algebras generated by two generic
$4\times 4$ matrices $X$ and $Y$. We
replace the variables $t_1,t_2$ with $x,y$ and denote $e_2=xy$. The Hilbert series of
$C$ and $T$ found by Berele and Stembridge \cite{BS} are
\[
h_C=\frac{P_C(x,y)}{Q_C(x,y)},\quad
h_T=\frac{P_T(x,y)}{Q_T(x,y)},\eqno(3)
\]
\[
P_C(x,y)=(1-e_2+e_2^2)(1-(x+y)e_2+(x+y)e_2^2+(x+y)^2e_2^2+(x+y)e_2^3-(x+y)e_2^4+e_2^6)
\]
\[
=(1-e_2+e_2^2)((1+e_2^3)^2-(x+y)e_2(1-e_2)^2(1+e_2)+(x^2+y^2)e_2^2),
\]
\[
Q_C(x,y)=(1-x)(1-x^2)(1-x^3)(1-x^4)(1-y)(1-y^2)(1-y^3)(1-y^4)
\]
\[
(1-xy)^2(1-x^2y)^2(1-xy^2)^2(1-x^3y)(1-xy^3)(1-x^2y^2)
\]
\[
=(1-xy)^3(1+xy)(1-x)^4(1+x)^2(1+x+x^2)(1+x^2)(1-e_2x)^2(1-e_2x^2)
\]
\[
(1-y)^4(1+y)^2(1+y+y^2)(1+y^2)(1-e_2y)^2(1-e_2y^2),
\]
\[
P_T(x,y)=1+e_2^2+e_2^3+e_2^5+(x+y)e_2^2,
\]
\[
Q_T(x,y)=(1-x)^2(1-x^2)(1-x^3)(1-y)^2(1-y^2)(1-y^3)
\]
\[
(1-xy)^2(1-x^2y)^2(1-xy^2)^2(1-x^3y)(1-xy^3)(1-x^2y^2)
\]
\[
=(1-xy)^3(1+xy)(1-x)^4(1+x)(1+x+x^2)(1-e_2x)^2(1-e_2x^2)
\]
\[
(1-y)^4(1+y)(1+y+y^2)(1-e_2y)^2(1-e_2y^2),
\]

The next two theorems give the multiplicity series of
these Hilbert series.

\begin{theorem}
The multiplicity series $m_C$
of the Hilbert series $H(C,x,y)$ of the
pure trace algebra of two generic $4\times 4$ matrices is
\[
m_C=
\frac{\alpha_4}{(1-t)^4}+\frac{\alpha_3}{(1-t)^3}+\frac{\alpha_2}{(1-t)^2}+\frac{\alpha_1}{1-t}
\]
\[
+ \frac{\beta_2}{(1+t)^2}+\frac{\beta_1}{1+t}+\frac{\gamma_0+\gamma_1t}{1+t+t^2}
+\frac{\delta_0+\delta_1t}{1+t^2}+\frac{\varepsilon_2}{(1-vt)^2}+\frac{\varepsilon_1}{1-vt}
+\frac{\varphi_0+\varphi_1t}{1-vt^2},
\eqno(4)
\]
where
\[
\alpha_4=\frac{(1-v+v^2)(1-v+v^2+4v^3+v^4-v^5+v^6)}{24(1-v)^{12}(1+v)^5(1+v+v^2)^2(1+v^2)},
\]
\[
\alpha_3=\frac{(1+v+v^2)(1-v+v^2)\alpha_3'}
{24(1-v)^{13}(1+v)^6(1+v+v^2)^4(1+v^2)^2},
\]
\[
\alpha_3'=3-4v-v^2+4v^3-3v^4-20v^5-12v^6-12v^7+7v^8+12v^9+v^{10}-4v^{11}+5v^{12},
\]
\[
\alpha_2=\frac{(1-v+v^2)\alpha_2'}{288(1-v)^{14}(1+v)^7(1+v+v^2)^4(1+v^2)^3},
\]
\[
\alpha_2'=59-97v-26v^2+223v^3+675v^4+840v^5+2501v^6+4049v^7+6799v^8+7754v^9
\]
\[
+6367v^{10}+3473v^{11}+2189v^{12}+768v^{13}+747v^{14}+271v^{15}-26v^{16}-97v^{17}+107v^{18},
\]
\[
\alpha_1=\frac{(1-v+v^2)\alpha_1'}{144(1-v)^{15}(1+v)^8(1+v+v^2)^5(1+v^2)^4}
\]
\[
\alpha_1'=34-86v-62v^2+106v^3+459v^4-624v^5-1887v^6-6630v^7-12804v^8-24712v^9
\]
\[
-40531v^{10}-57622v^{11}-62642v^{12}-57622v^{13}-40531v^{14}-24712v^{15}
-12804v^{16}\]
\[
-6630v^{17}-1887v^{18}
-624v^{19}+459v^{20}+106v^{21}-62v^{22}-86v^{23}+34v^{24},
\]
\[
\beta_2=\frac{1+v^4}{32(1-v)^6(1+v)^7(1+v^2)^3},
\]
\[
\beta_1=\frac{2-2v-4v^3+v^4-4v^5-2v^6-4v^7+v^8-4v^9-2v^{11}+2v^{12}}
{16(1-v)^7(1+v)^8(1+v+v^2)(1+v^2)^4(1-v+v^2)},
\]
\[
\gamma_0=\frac{1+2v}{9(1-v)^4(1+v)(1+v+v^2)^5(1-v+v^2)},
\]
\[
\gamma_1=\frac{1}{9(1-v)^3(1+v+v^2)^5(1-v+v^2)},
\]
\[
\delta_0= \frac{1}{8(1-v)^4(1+v)^3(1+v^2)^4},
\]
\[
\delta_1=0,
\]
\[
\varepsilon_2= \frac{v^{11}}{(1-v)^{14}(1+v)^7(1+v+v^2)^2(1+v^2)^3(1+v+v^2+v^3+v^4)},
\]
\[
\varepsilon_1= \frac{-v^{11}\varepsilon_1'}
{(1-v)^{15}(1+v)^8(1+v+v^2)^3(1+v^2)^4(1+v+v^2+v^3+v^4)^2(1-v+v^2)},
\]
\[
\varepsilon_1'=11+14v+31v^2+40v^3+60v^4+60v^5+72v^6
\]
\[
+60v^7+60v^8+40v^9+31v^{10}+14v^{11}+11v^{12},
\]
\[
\varphi_0=\frac{v^5\varphi_0'}{(1-v)^{15}(1+v)^3(1+v+v^2)^5(1+v+v^2+v^3+v^4)^2},
\]
\[
\varphi_0'=1+v+8v^2+10v^3+14v^4+17v^5+26v^6+17v^7
+14v^8+10v^9+8v^{10}+v^{11}+v^{12},
\]
\[
\varphi_1=\frac{2v^6(1-v+v^2)(1+2v^2+2v^3+2v^4+v^6)}
{(1-v)^{15}(1+v+v^2)^5(1+v+v^2+v^3+v^4)^2}.
\]
\end{theorem}

\begin{proof}
Direct calculations, which we have performed using Maple, show
that, replacing the expression (4) of $m_C$ instead of $M'(f)$
in the right hand side of (1) we obtain $h_C$ from (3). This completes
the proof, because of the injectivity of $M'$.
We want to say a couple of words how we have calculated $m_C$.
Applying Theorem 1 we know that $m_C$ is of the form
\[
m_C(t,v)=\frac{r(t,v)}{q(t,v)},
\]
where $r(t,v)\in {\mathbb C}(v)[t]$ and
\[
q(t,v)=(1-t)^4(1+t)^2(1+t+t^2)(1+t^2)(1-vt)^2(1-vt^2).
\]
Theorem 1 gives also that $\text{\rm deg}_tr\leq \text{\rm deg}_tq-2=12$.
Decomposing $m_C$ as a sum of elementary fractions with coefficients from ${\mathbb C}(v)$ and
denominators which are powers of $1-t$, $1+t$, $1+t+t^2$, $1+t^2$, $1-vt$, and $1-vt^2$, we
obtain that $m_C$ is as in (4). In order to calculate the coefficients
$\alpha_i,\beta_i,\gamma_i,\delta_i,\varepsilon_i,\varphi_i$, we replace
$h_C$ and $m_C$ in the left and right hind sides of (1), respectively.
Replacing $y$ with $e_2/x$, we apply Lemma 2 for
$K={\mathbb C}(e_2)$ and $\xi$ a zero of the denominator $q(x,e_2)$. In this way we obtain the
coefficients of the denominators of highest degree of
$1-t$, $1+t$, and $1-vt$. Replacing the obtained coefficients in (4), we calculate step by step
all the coefficients $\alpha_i,\beta_i,\varepsilon_i$. The coefficients
$\gamma_i,\delta_i,\varphi_i$, $i=0,1$, are obtained using both the zeros of $1+t+t^2$, $1+t^2$
and $1-vt^2$.
\end{proof}

The proof of the following theorem is similar.

\begin{theorem}
The multiplicity series $m_T$
of the Hilbert series $H(T,x,y)$ of the
mixed trace algebra of two generic $4\times 4$ matrices is
\[
m_T=
\frac{\alpha_4}{(1-t)^4}+\frac{\alpha_3}{(1-t)^3}+\frac{\alpha_2}{(1-t)^2}+\frac{\alpha_1}{1-t}
\]
\[
+\frac{\beta}{1+t}+\frac{\gamma_0+\gamma_1t}{1+t+t^2}
+\frac{\varepsilon_2}{(1-vt)^2}+\frac{\varepsilon_1}{1-vt}
+\frac{\varphi_0+\varphi_1t}{1-vt^2},
\]
where
\[
\alpha_4=\frac{1-v+3v^2-v^3+v^4}{6(1-v)^{12}(1+v)^3(1+v+v^2)^2},
\]
\[
\alpha_3=\frac{3-8v+4v^2-9v^3-8v^4-5v^5+12v^6-8v^7+7v^8}
{12(1-v)^{13}(1+v)^4(1+v+v^2)^3},
\]
\[
\alpha_2=\frac{\alpha_2'}{72(1-v)^{14}(1+v)^5(1+v+v^2)^4},
\]
\[
\alpha_2'=17-55v+124v^2+304v^3+540v^4+777v^5+1332v^6
\]
\[
+687v^7+468v^8+280v^9
+124v^{10}-73v^{11}+47v^{12},
\]
\[
\alpha_1=\frac{\alpha_1'}{144(1-v)^{15}(1+v)^6(1+v+v^2)^5}
\]
\[
\alpha_1'=25-134v+165v^2+123v^3-1758v^4-6240v^5-9439v^6-16537v^7-20250v^8
\]
\[
-16537v^9-9439v^{10}-6240v^{11}-1758v^{12}+123v^{13}+165v^{14}-134v^{15}
+25v^{16},
\]
\[
\beta=\frac{1-v+v^2-v^3+v^4}{16(1-v)^5(1+v)^6(1+v+v^2)(1+v^2)^2(1-v+v^2)},
\]
\[
\gamma_0=\frac{1+2v}{9(1-v)^4(1+v)(1+v+v^2)^5(1-v+v^2)},
\]
\[
\gamma_1=\frac{1}{9(1-v)^3(1+v+v^2)^5(1-v+v^2)},
\]
\[
\varepsilon_2= \frac{-v^8}{(1-v)^{14}(1+v)^5(1+v+v^2)^2(1+v^2)(1+v+v^2+v^3+v^4)},
\]
\[
\varepsilon_1= \frac{-v^8\varepsilon_1'}
{(1-v)^{15}(1+v)^6(1+v+v^2)^3(1+v^2)^2(1+v+v^2+v^3+v^4)^2(1-v+v^2)},
\]
\[
\varepsilon_1'=(3+3v+4v^2+4v^3+4v^4+3v^5+3v^6)
(3+v+5v^2+3v^3+5v^4+v^5+3v^6),
\]
\[
\varphi_0=\frac{2v^4\varphi_0'}{(1-v)^{15}(1+v)(1+v+v^2)^5(1+v+v^2+v^3+v^4)^2},
\]
\[
\varphi_0'=(2+2v+3v^2+3v^3+4v^4+v^5+v^6)
(1+v+4v^2+3v^3+3v^4+2v^5+2v^6)),
\]
\[
\varphi_1=\frac{v^4\varphi_1'}
{(1-v)^{15}(1+v+v^2)^5(1+v+v^2+v^3+v^4)^2},
\]
\[
\varphi_1'=1+5v+12v^2+18v^3+34v^4+37v^5+42v^6+37v^7+34v^8+18v^9+12v^{10}+5v^{11}+v^{12}.
\]
\end{theorem}

Theorems 3 and 4 give the explicit (but very complicated) form of the multiplicities
$m_{\lambda}(C)$ and $m_{\lambda}(T)$
for any $\lambda=(\lambda_1,\lambda_2)$.
We prefer to present the results in more compressed, but more convenient form.

\begin{proposition} The multiplicity series $m_C$ and $m_T$ are linear combinations of fractions
\[
\frac{v^at^b}{\pi^k(v)\rho^l(t,v)},
\eqno(5)
\]
where $0\leq a<\text{\rm deg}_v\pi(v)$,
$0\leq b<\text{\rm deg}_t\rho(t,v)$, and $\pi(v)$ and $\rho(t,v)$
are the following polynomials
\[
\pi(v)=1\pm v,1\pm v+v^2, 1+v^2, 1+v+v^2+v^3+v^4,
\]
\[
\rho(t,v)=1\pm t,1+t+t^2,1+t^2,1-vt,1-vt^2.
\]
The degrees of the denominators satisfy the inequality $k+l\leq 16$. The linear combinations
$M_C$ and $M_T$ of
fractions with denominators of degree $16$ of $m_C$ and $m_T$, respectively, are
\[
M_C=\frac{1}{(1-v)^{12}}\left(\frac{1}{2^83^2(1-t)^4}-\frac{1}{2^83^3(1-v)(1-t)^3}\right.
\]
\[
\left. +\frac{127}{2^{10}3^4(1-v)^2(1-t)^2}
-\frac{305}{2^93^5(1-v)^3(1-t)}\right.
\]
\[
\left.
-\frac{1}{2^{10}3^25(1-v)^2(1-vt)^2}-\frac{7}{2^93\cdot 5^2(1-v)^3(1-vt)}
+\frac{2^4(1+t)}{3^55^2(1-v)^3(1-vt^2)}\right),
\]
\[
M_T=16M_C.
\]
\end{proposition}

\begin{proof} Decomposing the rational functions
$\alpha_i,\beta_i,\gamma_i,\delta_i,\varepsilon_i,\varphi_i$
from Theorem 3 as linear combinations of
elementary fractions, we obtain that $m_C$ is a linear combination of fractions
of the form (5). We also see that the maximum degree $k+l$ is equal to 16
and this maximum is reached for the fractions
\[
\frac{1}{(1-v)^{16-l}(1-t)^l},\quad l=1,2,3,4,
\]
\[
\frac{1}{(1-v)^{16-l}(1-vt)^l}, \quad l=1,2,\quad
\frac{t^b}{(1-v)^{15}(1-vt^2)},\quad b=0,1.
\]
The explicit coefficients are found using Lemma 2.
The calculations for $m_T$ are similar, applying Theorem 4.
\end{proof}

The following theorem is the main result of our paper.
It is in the spirit of the description of
the multiplicities of the pure trace algebra of two $3\times 3$ matrices
given by Berele \cite{B1}.

\begin{theorem} Let $\lambda=(\lambda_1,\lambda_2)$ be
any partition.

{\rm (i)} The multiplicities $m_{\lambda}(C)$ of the pure trace
cocharacter of $4\times 4$ matrices satisfy the condition
\[
m_{\lambda}(C)=\begin{cases}
m_1+{\mathcal O}((\lambda_1+\lambda_2)^{13}), & \text{if $\lambda_1>3\lambda_2$,}\\
m_1+m_2+{\mathcal O}((\lambda_1+\lambda_2)^{13}), & \text{if $3\lambda_2\geq\lambda_1>2\lambda_2$,}\\
m_1+m_2+m_3+{\mathcal O}((\lambda_1+\lambda_2)^{13}), & \text{if $2\lambda_2\geq \lambda_1$,}\\
\end{cases}
\]
where
\[
m_1=\frac{(\lambda_1-\lambda_2)^3\lambda_2^{11}}{2^83^211!3!}
-\frac{(\lambda_1-\lambda_2)^2\lambda_2^{12}}{2^83^312!2!}
+\frac{127(\lambda_1-\lambda_2)\lambda_2^{13}}{2^{10}3^413!}
-\frac{305\lambda_2^{14}}{2^93^514!},
\]
\[
m_2=
\frac{(3\lambda_2-\lambda_1)^{14}}{2^{10}3^55^214!},
\]
\[
m_3=-\frac{(\lambda_1-\lambda_2)(2\lambda_2-\lambda_1)^{13}}{2^{10}3^25\cdot 13!}
-\frac{7(2\lambda_2-\lambda_1)^{14}}{2^93\cdot 5^214!}.
\]

{\rm (ii)} The multiplicities $m_{\lambda}(T)$ of the mixed trace
cocharacter of $4\times 4$ matrices satisfy
\[
m_{\lambda}(T)=16m_{\lambda}(C)+{\mathcal O}((\lambda_1+\lambda_2)^{13}).
\]
\end{theorem}

\begin{proof}
Let $\xi_1,\ldots,\xi_c\in {\mathbb C}$ be the zeros of $\pi(v)$ from (5).
Clearly, $\vert\xi_j\vert=1$. Then $1/\pi^k(v)$ is a linear combination of
$1/(1-\xi_jv)^s$, $s\leq k$. Using the formula
\[
\frac{1}{(1-\xi v)^s}=\sum_{q\geq 0}\binom{q+s-1}{s-1}\xi^qv^q,
\eqno(6)
\]
we obtain that
\[
\frac{1}{\pi^k(v)}=\sum_{q\geq 0}a_qv^q,
\]
where $\vert a_q\vert$ is bounded by a polynomial of degree $k-1$ in $q$.
A similar fact holds for $1/\rho^l(t,v)$ when $\rho=1\pm t,1+t+t^2, 1+t^2$.
Finally, the coefficients of the expansion
\[
\frac{1}{(1-vt)^2}=\sum_{p\geq 0}(p+1)(vt)^p
\]
are linear functions in $p$ and those of
\[
\frac{1}{1-vt}=\sum_{p\geq 0}(vt)^p,\quad
\frac{1}{1-vt^2}=\sum_{p\geq 0}(vt^2)^p
\]
are constants. In this way, the coefficients $b_{pq}$ of the expansion
\[
\frac{1}{\pi^k(v)\rho^l(t,v)}=\sum_{p,q\geq 0}b_{pq}t^pv^q
\]
are bounded by polynomials of degree $k+l-2$ in $p,q$ and satisfy
$b_{pq}={\mathcal O}((p+q)^{k+l-2})$. Since $k+l\leq 16$, the contribution
of maximum degree 14 to the multiplicities $m_{\lambda}(C)$ comes from the
expansion of $M_C$. Using (6) we obtain that $M_C$ is equal to
\[
\sum_{p,q\geq 0}\left(q^{11}\left(\frac{p^3}{2^83^211!3!}
-\frac{p^2q}{2^83^312!2!}+\frac{127pq^2}{2^{10}3^413!}
-\frac{305q^3}{2^93^514!}\right)+{\mathcal O}((p+q)^{13})\right)t^pv^q
\]
\[
+2^4(1+t)\sum_{r,s\geq 0}\frac{(s^{14}+{\mathcal O}((r+s)^{13}))(t^2v)^rv^s}{3^55^214!}
\]
\[
-\sum_{p,w\geq 0}\left(w^{13}\left(\frac{p}{2^{10}3^25\cdot 13!}
+\frac{7w}{2^93\cdot 5^214!}\right)
+{\mathcal O}((p+w)^{13})\right)(tv)^pv^w.
\]
Comparing the coefficients of the expansion of $m_C$ with those of $M_C$ we
obtain that
\[
m_{\lambda}(C)=q^{11}\left(\frac{p^3}{2^83^211!3!}
-\frac{p^2q}{2^83^312!2!}
+\frac{127pq^2}{2^{10}3^413!}
-\frac{305q^3}{2^93^514!}\right)
\]
\[
+\frac{2^4s^{14}}{3^55^214!}
-w^{13}\left(\frac{p}{2^{10}3^25\cdot 13!}
+\frac{7w}{2^93\cdot 5^214!}\right)+
+{\mathcal O}((\lambda_1+\lambda_2)^{13}),
\]
where $p,q,r,s,w\geq 0$ are such that
\[
\lambda_1-\lambda_2=p,\quad \lambda_2=q=r+s=p+w,
\]
\[
\lambda_1-\lambda_2=2s \text{ \rm or } \lambda_1-\lambda_2=2s+1.
\]
This gives the conditions
\[
s=\frac{1}{2}(\lambda_1-\lambda_2),\quad r=\frac{1}{2}(3\lambda_2-\lambda_1),
\]
if $\lambda_1-\lambda_2$ is even,
\[
s=\frac{1}{2}(\lambda_1-\lambda_2-1),\quad r=\frac{1}{2}(3\lambda_2-\lambda_1+1),
\]
if $\lambda_1-\lambda_2$ is odd,
and $w=2\lambda_2-\lambda_1$. The expression depending on $p,q$ gives the contribution $m_1$
for any $\lambda_1\geq \lambda_2$. The part with $r,s\geq 0$ gives $m_2$ when
$3\lambda_2\geq\lambda_1$. Finally, $p,w\geq 0$ gives $m_3$ when $2\lambda_2\geq\lambda_1$.
\end{proof}

\end{document}